\documentclass{amsart}

\usepackage{amscd}
\usepackage[all]{xy}
\usepackage{slashed}
\usepackage{todonotes}

\usepackage{comment}

\usepackage{hyperref}

\newtheorem{theo}{Theorem}[section]

\newtheorem{ques}{Question}
\newtheorem{conj}[theo]{Conjecture}

\theoremstyle{remark}
\newtheorem{rem}{Remark}[section]

\newcommand{\Z}{\mathbb{Z}}
\newcommand{\N}{\mathbb{N}}
\newcommand{\R}{\mathbb{R}}
\newcommand{\C}{\mathbb{C}}
\newcommand{\Q}{\mathbb{Q}}
\newcommand{\CP}{\mathbb{CP}}

\newcommand{\calM}{{\mathcal M}}

\newcommand{\calS}{{\mathcal S}}

\newcommand{\fraks}{\mathfrak{s}}

\newcommand{\Int}{\mathrm{Int}}
\newcommand{\Pin}{\mathrm{Pin}}
\newcommand{\Emb}{\mathrm{Emb}}
\newcommand{\Diff}{\mathrm{Diff}}
\newcommand{\BDiff}{B\mathrm{Diff}}
\newcommand{\EDiff}{E\mathrm{Diff}}

\newcommand{\Homeo}{\mathrm{Homeo}}
\newcommand{\BHomeo}{B\mathrm{Homeo}}
\newcommand{\Symp}{\mathrm{Symp}}
\newcommand{\TSymp}{\mathrm{TSymp}}

\newcommand{\del}{\partial}
\newcommand{\Ker}{\mathop{\mathrm{Ker}}\nolimits}

\newcommand{\id}{\mathrm{id}}
\newcommand{\SW}{\mathrm{SW}}
\newcommand{\SWbb}{\mathbb{SW}}

\newcommand{\Spinc}{\mathrm{Spin}^{c}}
\newcommand{\SO}{\mathrm{SO}}

\numberwithin{equation}{section}

\begin{document}

\title{Diffeomorphism groups and gauge theory for families}

\author{Hokuto Konno}
\curraddr{Graduate School of Mathematical Sciences, the University of Tokyo, 3-8-1 Komaba, Meguro, Tokyo 153-8914, Japan}
\email{konno@ms.u-tokyo.ac.jp}
\thanks{We would like to thank the referees for their valuable comments. The author was supported in part by JSPS KAKENHI Grant Numbers 25K00908.}

\dedicatory{Dedicated to Mikio Furuta}


\maketitle

\begin{abstract}
This article provides a survey of gauge theory for families, with a particular focus on its applications to diffeomorphism groups of $4$-manifolds that were developed during the period 2021--2025.
\end{abstract}

\section{Background}

We begin by explaining the background and motivation for our aims, without assuming any particular prior knowledge.

\subsection{Four-dimensional topology and gauge theory}

In the classification theory of manifolds, dimension four is special---
the establishment of this now well-known fact is regarded as one of the major achievements in the study of topology in the last century.

For example, many compact $4$-dimensional topological manifolds are known that admit infinitely many distinct smooth structures.
Here there are two relevant comparisons:
\begin{enumerate}
\item One is a comparison between dimension four and other dimensions. In dimensions other than four, it is known that a compact topological manifold admits at most finitely many smooth structures.

\item The other is a comparison between the topological category and the smooth category. The above fact shows that, in dimension four, there is a striking discrepancy between these two categories.
\end{enumerate}
These two comparisons are fundamental in four-dimensional topology.

The above fact concerning the number of smooth structures succinctly illustrates that the smooth category in dimension four is special.
In general, the smooth category in dimension four tends to exhibit distinctive forms of ``complexity''.
Such complexity typically manifests itself in the form that some quantity becomes infinite, or that some structure is unstable.

In order to prove such peculiarities of the smooth category in dimension four, one needs tools that can precisely capture smooth structures on $4$-manifolds.
Since Donaldson's work in the early 1980s, it has been well known that gauge theory provides such tools.
One of the typical arguments proceeds as follows, and is characterized by its essential use of nonlinear analysis:
first, one considers a certain nonlinear partial differential equation (originating from gauge theory in physics) on a smooth $4$-manifold.
The space of solutions admits certain symmetries arising from the equation (gauge symmetries).
By suitably ``counting'' the moduli space obtained as the quotient of the solution space by these symmetries, one obtains numerical invariants of the $4$-manifold.
These invariants are then used to distinguish $4$-manifolds that are homeomorphic but not diffeomorphic.

However, to carry out this idea, it is necessary to compute gauge-theoretic invariants.
This is often not easy, and in practice explicit computations are possible only in limited situations.
More than forty years have passed since Donaldson's work, but it seems that we have captured only a small portion of the peculiarities of the smooth category in dimension four.

\subsection{Diffeomorphism groups of higher-dimensional manifolds}

Looking back to the period before the emergence of gauge theory, the classification theory of higher-dimensional (five-dimensional and above) manifolds reached a certain milestone in the late 1960s with the work of Kirby and Siebenmann~\cite{KS77}.
For example, the finiteness of smooth structures in dimensions other than four mentioned above is (in higher dimensions) a consequence of their work.
At the heart of their method lies, roughly speaking, a kind of surgical technique called Whitney's trick.
This provides a bridge that reduces geometric problems about manifolds to algebraic problems.
Once reduced to algebra, geometric problems in higher dimensions become manageable to a certain extent.
Whitney's trick is an argument that removes intersections in a manifold along $2$-dimensional disks, and it is the main reason why dimension four is excluded in the classification theory of higher-dimensional manifolds (the number four appears because, in manifolds of dimension at least five, one can eliminate self-intersections of disks by a general position argument, whereas this is impossible in dimension $4=2+2$).

On the other hand, in manifold topology, the automorphism group of a manifold---the diffeomorphism group---has always been a fundamental object of interest. For example, in the case of two-dimensional manifolds, the study of mapping class groups has long occupied a central place in low-dimensional topology. In contrast, since the beginning of this century, the study of diffeomorphism groups of higher-dimensional manifolds has made remarkable progress, some of which we will discuss later in this survey. The methods involved combine sophisticated techniques from algebraic topology and differential topology, and the differential-topological techniques again include Whitney's trick. Recent progress in the study of diffeomorphism groups of higher-dimensional manifolds recalls the golden age of topology in the 1960s.

\subsection{Diffeomorphism groups of four-dimensional manifolds}

Let us return to dimension four.
Compared with the flourishing developments in other dimensions, knowledge about diffeomorphism groups of $4$-manifolds had been extremely limited until recently.
This stems from the absence of Whitney's trick, for reasons similar to those that make the classification theory of $4$-manifolds difficult.

However, in the classification theory of $4$-manifolds, there exist phenomena that are not merely difficult, but genuinely different from those in other dimensions.
It is therefore natural to ask whether the same is true for diffeomorphism groups.
In addition to comparisons with other dimensions, it is also natural to study diffeomorphism groups from another fundamental viewpoint in four-dimensional topology, namely the comparison between the topological and smooth categories---that is, to compare homeomorphism groups and diffeomorphism groups:

\begin{ques}
Investigate whether diffeomorphism groups of $4$-manifolds possess properties that differ from those in other dimensions.
\end{ques}

\begin{ques}
\label{ques: Diff vs Homeo}
Investigate the difference between the homeomorphism group and the diffeomorphism group of a $4$-manifold.
\end{ques}

When considering these questions, the analogy with the classification theory of manifolds suggests that one should expect that diffeomorphism groups of $4$-manifolds exhibit singularly complex behavior.
The first issue is how to formulate ``complexity'' for diffeomorphism groups.
Compared with the classification theory of manifolds, complexity in diffeomorphism groups manifests itself in a wide variety of ways.
In this article, we formulate complexity in several different ways and explain that complexities specific to the smooth category in dimension four do indeed appear in diffeomorphism groups.

\subsection{Approach: gauge theory for families}

To address the above problems, gauge theory for families, which develops gauge theory for smooth families of $4$-manifolds, provides an effective approach.
A typical argument proceeds as follows.
Suppose we are given a smooth family of $4$-manifolds, that is, a smooth fiber bundle whose fibers are $4$-manifolds.
On each fiber, a gauge-theoretic partial differential equation is defined.
One considers the moduli spaces (quotients of solution spaces) assembled fiberwise over the base space.
This assembled moduli space is called the {\it parameterized moduli space}.
By ``counting'' this parameterized moduli space, one obtains invariants of the family of $4$-manifolds.
Using these invariants to distinguish fiber bundles that are typically isomorphic as topological bundles but not as smooth bundles, one can obtain ``lower bounds'' on the complexity of the diffeomorphism group of the fiber, which serves as the structure group of the smooth fiber bundle.

The first application of this idea to diffeomorphism groups of $4$-manifolds was due to Ruberman~\cite{Rub98} in 1998.
In this pioneering work, he implemented the idea in the case where the base space of the family is one-dimensional, and provided one answer to Question~\ref{ques: Diff vs Homeo} (see Section~\ref{section exotic Dehn twist}) by constructing the first example of an exotic diffeomorphism of a $4$-manifold, namely, a diffeomorphism that is topologically isotopic but not smoothly isotopic to the identity.
Over the subsequent decade, there were follow-up works by Ruberman himself~\cite{Rub99,Rub01} and related works by Nakamura~\cite{Naka03,Naka10}, but applications of gauge theory for families to diffeomorphism groups remained sporadic.

Since the mid-2010s, a more systematic study of gauge theory for families has been pursued, and recent years have seen rapid developments by many researchers.
As a result, new perspectives on the above questions are now emerging, and this article surveys part of that progress.

We briefly summarize the contents of this article.
In Section~\ref{section homological instavility}, we explain the homological instability of moduli spaces of $4$-manifolds.
In Section~\ref{section infiniteness}, we discuss various kinds of infiniteness associated with diffeomorphism groups of $4$-manifolds.
In Section~\ref{section exotic Dehn twist}, we summarize recent results on exotic diffeomorphisms of $4$-manifolds.
In Section~\ref{section others}, we briefly survey other topics related to gauge theory for families, as well as another important recent development, namely Kontsevich characteristic classes.
Finally, in Section~\ref{section future}, we touch very briefly on future problems.

Applications of gauge theory for families to diffeomorphism groups up to around March 2022 are summarized in the survey~\cite{KonnoRonsetsu}.
In this article, we restrict attention to results obtained after that point and up to around April 2025.

\section{Homological instability}
\label{section homological instavility}

\subsection{Moduli spaces of manifolds}

We begin by recalling the classifying space of the diffeomorphism group.
For a compact smooth manifold $X$, let $\Diff(X)$ denote the diffeomorphism group of $X$.
The group $\Diff(X)$ carries a natural topology, namely the $C^\infty$ topology, and becomes a topological group under composition.
Its classifying space $\BDiff(X)$, intuitively, parametrizes all manifolds diffeomorphic to $X$, and is called the moduli space of manifolds diffeomorphic to $X$.

It is also easy to describe an explicit model for $\BDiff(X)$.
If we write $\R^\infty$ for the direct limit of the natural inclusions
\[
\R^{1} \subset \R^{2} \subset \cdots,
\]
then one has
\[
\BDiff(X) = \Emb(X,\R^\infty)/\Diff(X).
\]
Here $\Emb(X,\R^\infty)$ denotes the space of smooth embeddings of $X$ into $\R^\infty$, and $\Diff(X)$ acts on $\Emb(X,\R^\infty)$ by reparametrization.
An embedding of $X$ into $\R^\infty$ may be regarded as an embedding into a Euclidean space $\R^N$ of sufficiently large dimension, and ``reparametrization'' means precomposing an embedding $X \hookrightarrow \R^\infty$ with a diffeomorphism of its domain.
In intuitive terms, $\BDiff(X)$ is the space obtained by collecting all submanifolds of $\R^\infty$ that are diffeomorphic to $X$.
This viewpoint parallels the description of the classifying space of the orthogonal group $O(n)$ as the infinite-dimensional Grassmannian of $n$-dimensional linear subspaces of $\R^\infty$.

The reason $\BDiff(X)$ is called the moduli space is that it classifies families of $X$, that is, smooth fiber bundles with fiber $X$.
Indeed, by the general theory of classifying spaces, for any sufficiently nice topological space $B$ (for example, a CW complex), there is a one-to-one correspondence
\[
[B, \BDiff(X)] \overset{1:1}{\leftrightarrow} \left\{\text{isomorphism classes of fiber bundles over $B$ with fiber $X$}\right\}.
\]
When $X$ has boundary, one instead considers the group $\Diff_\del(X)$ of diffeomorphisms that are the identity near the boundary of $X$, together with its classifying space.

Just as, in the classification of vector bundles, it is important to understand the topology of Grassmann manifolds, the topology of the moduli space $\BDiff(X)$ is among the most basic objects of study in the topology and geometry of fiber bundles.
Among its topological features, the cohomology of the moduli space corresponds to characteristic classes of fiber bundles with fiber $X$, and is one of the central objects of interest:
\[
H^\ast(\BDiff(X)) \overset{1:1}{\leftrightarrow} \left\{\text{characteristic classes of fiber bundles with fiber $X$}\right\}.
\]

That the study of the topology of $\BDiff(X)$ has a scope quite different from the classification theory of manifolds is already evident from the history of two-dimensional topology.
As is well known, the classification of surfaces was completed more than a century ago.
However, the moduli space $\BDiff(X)$ of surfaces is naturally related to the moduli space of Riemann surfaces, and its topology is closely tied to the study of mapping class groups.
This remains an active area of research today and will continue to do so.
Similarly, for higher-dimensional manifolds, the topology of $\BDiff(X)$ has been studied in earnest since the beginning of this century, going well beyond the classification theory of manifolds.

It may be helpful to think of the study of the moduli space of manifolds as a ``higher'' analogue of the classification of manifolds.
Indeed, the object most deserving of the name ``the moduli space of manifolds'' is the disjoint union, taken over all diffeomorphism classes $[X]$ of smooth manifolds,
\[
\calM = \bigsqcup_{[X]} \BDiff(X).
\]
The classification theory of manifolds is the study of the set of connected components $H_{0}(\calM)$ or $\pi_{0}(\calM)$, whereas the study of $\BDiff(X)$ is, in many contexts, naturally viewed as the study of higher homology groups $H_{k}(\calM)$ or higher homotopy groups $\pi_{k}(\calM)$ of a given connected component, and related structures.
For example, the recent results of the author~\cite{konno2023homology} and of Auckly--Ruberman~\cite{auckly2025familiesdiffeomorphismsembeddingspositive}, discussed later, show that performing the operation of connected sum with $S^{2}\times S^{2}$ on a $4$-manifold $k$ times ($k>0$) gives rise to interesting elements of $H_{k}(\calM)$ or $\pi_{k}(\calM)$ from suitable collections of elements in $H_{0}(\calM)$ or $\pi_{0}(\calM)$.
The prototype of this argument is the pioneering work of Ruberman~\cite{Rub98}, which produces an element of $\pi_{1}(\calM)$ from several elements of $\pi_{0}(\calM)$.

\subsection{Homological stability in dimensions other than four}
\label{subsection: homological stability}

A fundamental problem is to compute the cohomology $H^\ast(\BDiff(X))$, or dually the homology $H_\ast(\BDiff(X))$, but this is in general extremely difficult.
For instance, even when $X$ is a surface, we are still far from a complete computation.
However, if one stabilizes a manifold in an appropriate way, then, remarkably, there are cases in which $H_\ast(\BDiff(X))$ becomes computable.
The key notion is that of {\it homological stability}, which we explain below.

Let $W$ be a compact smooth manifold of even dimension $2n$ with boundary.
We adopt the following model for the connected sum $W\# S^n\times S^n$: we glue $W$ and $S^n\times S^n \# \del W\times[0,1]$ (where $\#$ denotes the interior connected sum) by identifying $\del W$ with $\del W\times\{0\} \subset S^n\times S^n \# \del W\times[0,1]$.
Then, by extending diffeomorphisms by the identity, we obtain an injective homomorphism
\[
s : \Diff_\del(W) \hookrightarrow \Diff_\del(W\# S^n\times S^n).
\]
The map $s$ and the maps induced by it are called stabilization maps.
The following result is known as homological stability for moduli spaces of manifolds. It is due to Harer~\cite{Harer85} in the two-dimensional case and to Galatius--Randal-Williams~\cite{Galatius2018} in dimensions at least six.
The theorem of Galatius--Randal-Williams provides a foundation for the recent study of diffeomorphism groups of high-dimensional manifolds.
In what follows, for $N>0$, we write $\#_{N}S^n\times S^n$ for the connected sum of $N$ copies of $S^{n}\times S^{n}$.

\begin{theo}[Harer \cite{Harer85}, Galatius--Randal-Williams \cite{Galatius2018}]
\label{thm: homological stability}
Let $W$ be a compact smooth simply-connected manifold of even dimension $2n \neq 4$.
Fix $k\geq 0$.
Then, for any $N \gg k$, the map induced by $s$,
\[
s_\ast : H_k(\BDiff_\del(W\#_{N}S^n\times S^n);\Z)
\to 
H_k(\BDiff_\del(W\#_{N+1}S^n\times S^n);\Z),
\]
is an isomorphism.
\end{theo}

By this theorem, for sufficiently large $N$, $H_k(\BDiff_\del(W\#_N S^n\times S^n);\Z)$ can be identified with the direct limit
\[
\lim_{N \to +\infty}H_k(\BDiff_\del(W\#_N S^n\times S^n);\Z).
\]
This limit is called the {\it stable homology}, and it is an object that can be computed by homotopy-theoretic methods.
Over $\Q$, it can be described explicitly and is generated by the well-known characteristic classes called the Mumford--Morita--Miller classes.
The solution of Mumford's conjecture by Madsen--Weiss~\cite{MadsenWeiss07} establishes this in dimension two, and later Galatius--Randal-Williams~\cite{Galatius2018} extended it to higher dimensions.

The reason dimension four is excluded from the above theorem of Galatius--Randal-Williams is that Whitney's trick is used at a key point in the proof.
As already mentioned, the failure of Whitney's trick has been the main obstacle to extending higher-dimensional theorems on the classification of manifolds to dimension four.
Likewise, in the study of diffeomorphism groups, higher-dimensional arguments again fail to extend to dimension four precisely because Whitney's trick is unavailable there.

\subsection{Homological instability in dimension four}

As we saw in Theorem~\ref{thm: homological stability} and the subsequent discussion, in even dimensions other than $4$, after sufficiently many stabilizations, one can compute the homology $H_k(\BDiff_\del(W))$ of the moduli space.
The following theorem asserts that, in dimension four, the key step---namely, homological stability---fails.
Here we restrict attention to manifolds with boundary obtained by puncturing closed manifolds.
That is, in the theorem below, for a closed $4$-manifold $X$, we consider $\mathring{X} = X \setminus \Int(D^4)$.

\begin{theo}[Konno--Lin \cite{konno2022homological}]
\label{theo: homological instability}
Let $X$ be a smooth simply-connected closed $4$-manifold.
Fix $k>0$.
Then there exists a sequence $0<N_1 < N_2 < \cdots$ with $N_i \to +\infty$ such that, for each $N_i$,
\[
s_\ast : H_k(\BDiff_\del(\mathring{X}\#_{N_i}S^2 \times S^2);\Z)
\to 
H_k(\BDiff_\del(\mathring{X}\#_{N_i+1}S^2 \times S^2);\Z)
\]
is not an isomorphism.
\end{theo}

In other words, in contrast to the situation in dimensions other than four, the homology of the moduli space never stabilizes, regardless of how many times one takes the connected sum with $S^2\times S^2$.
Moreover, the same conclusion holds with ``isomorphism'' replaced by ``injective'' or ``surjective'' \cite{konno2022homological}: characteristic classes are created and annihilated infinitely many times.

In contrast to the situation in other dimensions described in Section~\ref{subsection: homological stability}, Theorem~\ref{theo: homological instability} suggests that, in dimension four, the computation of the homology of moduli spaces is difficult to reduce to homotopy-theoretic methods.
A phenomenon similar to that in the classification theory of manifolds occurs here as well.

In the topological category, the following is currently known.
In dimension two, there is no difference between the topological and smooth categories, and thus Harer~\cite{Harer85} already implies homological stability for the classifying space $\BHomeo_\del(W)$ of the homeomorphism group.
In higher dimensions, the corresponding stability result in the topological category was proved by Kupers~\cite{kupers2016provinghomologicalstabilityhomeomorphisms}.
In dimension four, homological stability in the topological category has not yet been established.

In the proof of Theorem~\ref{theo: homological instability}, one constructs explicit unstable elements in the homology of the moduli space, but all of them become trivial after passing to the topological category.
Thus, the proof also reveals a difference between the homology of moduli spaces in the smooth category and its counterpart in the topological category.

\begin{rem}
Nariman~\cite{Nariman17a,Nariman17} established a version of homological stability for $\Diff^{\delta}(X)$, the diffeomorphism group endowed with the discrete topology.
In dimension four, however, this discrete analogue also fails~\cite{konno2022homological}.
\end{rem}

\subsection{Proof of homological instability (1): constructing gauge-theoretic characteristic classes}

Homological instability (Theorem~\ref{theo: homological instability}) is proved by constructing and computing characteristic classes using gauge theory for families.
Let $X$ be an oriented smooth closed $4$-manifold, and let $k\geq 0$.
Let $\Diff^+(X)$ denote the group of orientation-preserving diffeomorphisms.
In \cite{konno2022homological}, a $\Z/2$-valued characteristic class
\begin{align}
\label{eq: SW char}
\SWbb^k(X) \in H^k(\BDiff^+(X);\Z/2)
\end{align}
was introduced under the assumption that $b^+(X)\geq k+2$.
This class is constructed by combining a gauge-theoretic characteristic class defined by the author \cite{K21} with an idea of Ruberman~\cite{Rub01} for a numerical invariant associated with a diffeomorphism
\footnote{\cite{K21} considers the diffeomorphism group $\Diff(X,\fraks)$ preserving a spin$^c$ structure $\fraks$, and defines classes $\SWbb^\ast(X,\fraks)\in H^\ast(\BDiff(X,\fraks))$ using families of the Seiberg--Witten equations.
To treat $\Diff^+(X)$ rather than $\Diff(X,\fraks)$, one uses a slight modification of the idea in \cite{Rub01}.
In \cite{Rub01}, for a diffeomorphism $f\in \Diff^+(X)$ that need not preserve $\fraks$, a numerical invariant $\SW(X,f,\fraks)\in \Z/2$ (and, for some $f$, $\SW(X,f,\fraks)\in \Z$) is defined.
The invariant $\SW(X,f,\fraks)$ is defined by considering the Seiberg--Witten equations for all spin$^c$ structures of the form $\{(f^{n})^\ast \fraks\}_{n\in \Z}$.
By contrast, the characteristic class defined in \cite{konno2022homological} is obtained by considering the Seiberg--Witten equations for all spin$^c$ structures of a fixed formal dimension.
The invariant $\SW(X,f,\fraks)$ is not a homomorphism with respect to composition in $f$, but by considering all spin$^c$ structures one gains a homomorphism property and obtains a cohomology class on $\BDiff^{+}(X)$.}.
Very roughly speaking, it is defined by counting solutions of families of the Seiberg--Witten equations
\begin{align*}
\left\{
\begin{array}{ll}
F^+_A = \sigma(\Phi, \Phi),\\
\slashed{D}_A \Phi = 0,
\end{array}
\right.
\end{align*}
over each $k$-cell of $\BDiff^{+}(X)$.

Let us say a little more about the idea behind the construction of the characteristic class $\SWbb^k(X)$.
The Seiberg--Witten equations can be written down once one fixes a topological datum, i.e., a spin$^c$ structure on $X$, together with a Riemannian metric.
Let us consider these data in families.
To discuss the most general setting, consider the universal $X$-bundle
\[
E := \EDiff^+(X)\times_{\Diff^+(X)}X \to \BDiff^+(X)
=:B.
\]
The monodromy of $E$ acts on the set $\Spinc(X)$ of isomorphism classes of spin$^c$ structures on $X$ by pullback, and hence one obtains an associated family over $B$ with fiber $\Spinc(X)$.
Writing this family as $\Spinc(E)\to B$, we see that $\Spinc(E)\to B$ is a covering space of $B$.
In general, the action of $\Diff^{+}(X)$ on $\Spinc(X)$ is complicated, and $\Spinc(E)\to B$ is a highly nontrivial covering space.

Choose a smooth family of Riemannian metrics on $E$.
Then, over each point $b\in B$, the fiber $E_b$ carries a metric $g_b$, together with the fiber $\Spinc(E_b)$ of $\Spinc(E)$.
For each element of $\Spinc(E_b)$, one can consider the Seiberg--Witten equations with respect to the metric $g_b$.
Take the disjoint union of these equations over all elements of $\Spinc(E_b)$.
As $b$ varies, we obtain a family over $B$ of systems of Seiberg--Witten equations, in total $\#\Spinc(X)$ such systems.
Although on each fiber we are merely considering a collection of unrelated Seiberg--Witten equations, over the entire base $B$ one obtains a family of equations that is glued together in an intricate way along the covering space $\Spinc(E)\to B$.

Since one can functorially associate to any topological space a weakly homotopy equivalent CW complex, we may assume without loss of generality that $B=\BDiff^+(X)$ is equipped with a CW complex structure.
Over each $k$-cell of $B$, we count the moduli spaces for the above family of systems of Seiberg--Witten equations, and thereby obtain a number.
In this way we construct a cochain on $B$, and one can show that it is in fact a cocycle.
Its cohomology class is the characteristic class \eqref{eq: SW char}, and one can show that it does not depend on the choice of the family of metrics (or on the perturbations needed for technical reasons).
We use $\Z/2$ coefficients because the monodromy of $E$ acts nontrivially on the orientation (homology orientation) of the moduli spaces.
More precisely, in the proof of Theorem~\ref{theo: homological instability}, one selects only spin$^c$ structures of an appropriate formal dimension, and, in order to obtain nontrivial computations, one uses a version of the above construction quotiented by a $\Z/2$-symmetry called charge conjugation.

An important property of the characteristic class $\SWbb^k(X)$ is that it is unstable.
That is, one can prove that it vanishes after applying the pullback $s^\ast$ of the stabilization map (more precisely, of a finite iterate of it).
This is a families analogue of the vanishing of Seiberg--Witten invariants under suitable connected sums, including connected sums with $S^2\times S^2$.
Thanks to this instability, once one can compute the characteristic class as discussed next, one sees that gauge-theoretic characteristic classes provide a tool for capturing the complexity of diffeomorphism groups of $4$-manifolds in the form of homological instability.

\subsection{Proof of homological instability (2): computing the characteristic classes}

To use the above characteristic classes $\SWbb^k(X)$, one needs to find fiber bundles for which $\SWbb^k(X)$ is nontrivial, and such families are constructed as follows.
This is a higher-dimensional analogue of the argument carried out by Ruberman~\cite{Rub98} for families over $S^1$.
More precisely, it is the Seiberg--Witten analogue in higher dimensions of an argument of Ruberman in \cite{Rub98}, which used Yang--Mills theory.

For a closed smooth oriented $4$-manifold $M$, consider a $4$-manifold of the form $M\#_{k}S^2\times S^2$.
Choose $k$ mutually commuting diffeomorphisms $f_1,\ldots,f_k$ appropriately; roughly speaking, each is supported in a copy of $S^2\times S^2$ and acts on its homology by multiplication by $-1$.
Then the multiple mapping torus of these diffeomorphisms gives a family of $M\#_{k}S^2\times S^2$ over $T^k$.

If we take $M$ to be a symplectic manifold satisfying suitable properties, then for this family over $T^k$ the characteristic class \eqref{eq: SW char} becomes nontrivial.%
\footnote{Briefly, the required property of $M$ is that the sum of Seiberg--Witten invariants over spin$^c$ structures other than the one coming from the symplectic structure, even after quotienting by charge conjugation, vanishes mod $2$.}
At the heart of the proof of nontriviality lies a computation formula for families Seiberg--Witten invariants for families of $4$-manifolds due to Baraglia and the author~\cite{BK20gluing}.
This is a Taubes-type gluing formula, and can be viewed as a families version of the gluing argument for solutions of the Seiberg--Witten equations that appears in the proof of blow-up formulas for Seiberg--Witten invariants; see, for instance, \cite{Ni20}.
In conclusion, our gluing formula reduces the computation of the families Seiberg--Witten invariants of the family of $M\#_{k}S^2\times S^2$ to the computation of the Seiberg--Witten invariants of $M$.

This gluing formula computes the contribution of parametrized moduli spaces associated with monodromy-invariant spin$^{c}$ structures to the characteristic class.
The contributions coming from individual spin$^{c}$ structures that are not monodromy-invariant are in general highly complicated and seem difficult to compute directly.
However, if one considers the orbits of the monodromy action on the set of all spin$^{c}$ structures, then one finds that the contributions from all orbits other than the monodromy-invariant ones cancel over $\Z/2$.
This is a consequence of having chosen simple diffeomorphisms on $S^{2}\times S^{2}$ that act on homology by multiplication by $-1$.

Putting everything together, one concludes that only the Seiberg--Witten invariants of $M$ ultimately affect the computation of the characteristic class, and hence one obtains the desired nontriviality.
The contribution of the Seiberg--Witten invariants of $M$ is given by $\SWbb^{0}(M)$.
Therefore, to ensure nontriviality of the characteristic class, one must choose $M$ so that a suitably defined sum of the Seiberg--Witten invariants of $M$ over spin$^{c}$ structures is nontrivial.
For this, one needs information more refined than Taubes' computation of the Seiberg--Witten invariant for the canonical spin$^{c}$ structure associated with a symplectic structure~\cite{Taubes94sympl}.

This completes the core of the proof of the homological instability theorem.
To finish the proof, it suffices, for a given smooth simply-connected closed $4$-manifold $X$, to find a sequence $N_i\to +\infty$ and a sequence of symplectic $4$-manifolds $M_i$ satisfying the above properties such that
\[
X\#_{N_i}S^2\times S^2 \cong M_i\#_k S^2\times S^2.
\]
This reduces to a certain kind of geography problem.
Once this is achieved, combining it with the instability of the characteristic class $\SWbb^k(X)$ noted above shows that the stabilization maps $s_{\ast}$ appearing in Theorem~\ref{theo: homological instability} are not surjective.

That the stabilization maps $s_{\ast}$ are not injective is seen as follows.
Consider the family $E$ constructed as above and another family $E'$ obtained by the same construction but with $M$ replaced by a $4$-manifold $M'$ whose Seiberg--Witten invariants are trivial.
By the same gluing argument as above, one finds that the characteristic class $\SWbb^k(X)$ is trivial for $E'$.
It follows that the elements of $H_{k}(\BDiff(M\#_{k} S^{2} \times S^{2});\Z)$ determined by the classifying maps of $E$ and $E'$ are distinct.
On the other hand, by choosing $M$ and $M'$ appropriately, one may assume in advance that $M\# S^{2} \times S^{2}$ and $M'\# S^{2} \times S^{2}$ are diffeomorphic, and then it is easy to see that $E$ and $E'$ become isomorphic fiber bundles after taking the fiberwise connected sum with the trivial bundle $(S^{2} \times S^{2}) \times T^{k}$ along a trivial $D^{4}$-subbundle.%
\footnote{From the construction, both $E$ and $E'$ contain the trivial bundle $D^{4} \times T^{k}$ as a subbundle.
Along this trivial bundle, one can take the fiberwise connected sum with the trivial bundle $(S^{2} \times S^{2}) \times T^{k}$.}
In other words, two distinct elements of $H_{k}(\BDiff(M\#_{k} S^{2} \times S^{2});\Z)$ become equal after applying $s_{\ast}$ and mapping to $H_{k}(\BDiff(M\#_{k+1} S^{2} \times S^{2});\Z)$.
This shows that $s_{\ast}$ is not injective.

\section{Infiniteness phenomena for diffeomorphism groups}
\label{section infiniteness}

\subsection{Finiteness in higher-dimensions}

In most cases, the homology and homotopy groups of the diffeomorphism group $\Diff(X)$ and of its classifying space are infinite.
Given an infinite group, one of the first structural questions to ask is whether it is finitely generated, which is a weak form of finiteness.
In dimensions other than four, strong finiteness theorems are known under suitable assumptions on the fundamental group.
For example, no counterexample to the following conjecture is currently known:

\begin{conj}
\label{conj: finiteness}
Let $X$ be a smooth closed manifold.
If $\pi_1(X)$ is finite and $\dim X\neq 4$, then $\pi_k(\Diff(X))$ and $H_k(\BDiff(X);\Z)$ should be finitely generated for all $k\geq 0$.
\end{conj}

Furthermore, in even dimensions this conjecture has already been established (Kupers~\cite{Kupers19}, Bustamante--Krannich--Kupers~\cite{bustamante2023finiteness}).
Their proofs use homological stability, namely the part of Theorem~\ref{thm: homological stability} valid in dimensions at least six.
Looking back, dimension four was excluded from the proof of homological stability because of the failure of Whitney's trick.
Thus, here as well, the failure of Whitney's trick makes it unclear whether such finiteness theorems hold in dimension four.

\subsection{Infiniteness of higher homotopy groups of $\Diff(X^{4})$}

It is natural to ask whether a four-dimensional analogue of Conjecture~\ref{conj: finiteness} fails.
As mentioned earlier, unlike in other dimensions, compact $4$-manifolds can admit infinitely many exotic smooth structures.
Searching for counterexamples to a four-dimensional analogue of Conjecture~\ref{conj: finiteness} can be viewed as asking whether there is an analogue, at the level of automorphism groups, of this phenomenon -- namely, a certain kind of infiniteness specific to dimension four.
For the homotopy groups of $\Diff(X)$, the following result is the first in this direction:

\begin{theo}[Baraglia \cite{B21}]
\label{theo: Baraglia}
$\pi_1(\Diff(K3))$ is not finitely generated.
\end{theo}

More strongly, \cite{B21} shows that $\pi_1(\Diff(K3))$ contains $\Z^{\oplus \N}$ as a direct summand.
This result is proved using families Seiberg--Witten invariants.
This result was later generalized and refined by Lin~\cite{lin2022family} to a broad class of $4$-manifolds, including all elliptic surfaces and complete intersections.

On the other hand, based on examples quite different from those of Baraglia and Lin, Auckly and Ruberman obtained results extending Ruberman's series of works on $1$-parameter families~\cite{Rub98,Rub99,Rub01} to higher homotopy groups.
In particular:

\begin{theo}[Auckly--Ruberman~\cite{auckly2025familiesdiffeomorphismsembeddingspositive}]
\label{theo: AucklyRuberman}
For each $k\geq 1$, there exists a smooth simply-connected closed $4$-manifold $X$ such that $\pi_k(\Diff(X))$ is not finitely generated.
\end{theo}

More strongly, \cite{auckly2025familiesdiffeomorphismsembeddingspositive} proved that, for each $k\geq 1$, there exists such an $X$ for which $\pi_k(\Diff(X))$ contains $\Z^{\oplus \N}$ as a direct summand.
A qualitative difference between the examples in \cite{B21,lin2022family} and those in \cite{auckly2025familiesdiffeomorphismsembeddingspositive} is that, for the former, it is not clear whether the families over $S^{2}$ considered in \cite{B21,lin2022family}, corresponding to elements of $\pi_{1}(\Diff(X))$, are topologically trivial.
At least, the explicit generators constructed in \cite{B21,lin2022family} are not topologically trivial, since the vertical tangent bundle has nontrivial rational Pontryagin classes, which are also defined in the topological category.
On the other hand, the families constructed in \cite{auckly2025familiesdiffeomorphismsembeddingspositive} are topologically trivial, and hence yield a $\Z^{\oplus \N}$-summand in
\[
\Ker(\pi_{k}(\Diff(X)) \to \pi_{k}(\Homeo(X))).
\]

Theorem~\ref{theo: AucklyRuberman} is again proved by computing families Seiberg--Witten invariants.
The idea is to extend inductively the computation of invariants for $1$-parameter families due to Ruberman~\cite{Rub98} (and its Seiberg--Witten counterpart~\cite{BK20gluing}) to the setting in which one takes connected sums with as many copies of $S^{2}\times S^{2}$ as dictated by the degree of the homotopy group under consideration, and then computes invariants over spheres of arbitrary dimension.

\subsection{Infiniteness of mapping class groups and the homology of $\BDiff(X^{4})$}

The results discussed above concern higher homotopy groups of diffeomorphism groups.
Next, we discuss infiniteness phenomena for the mapping class group of $X$ and for the homology of $\BDiff(X)$.
This will give another application of gauge-theoretic characteristic classes.
As background, the following classical result concerns finiteness of mapping class groups outside dimension four:

\begin{theo}[Sullivan \cite{Sull77}]
\label{theo: Sullivan}
Let $X$ be a smooth simply-connected closed manifold of $\dim \geq 5$.
Then $\pi_0(\Diff(X))$ is finitely generated.
\end{theo}

In fact, Sullivan proved a stronger finiteness statement; for example, under the same assumption $\pi_0(\Diff(X))$ is finitely presented \footnote{In \cite{Sull77} the case of dimension five is not stated explicitly, but it is known that the analogous statement also holds in dimension five.}.
Of course, the analogous statement is true in dimensions at most three as well, so for simply-connected manifolds one concludes that the mapping class group is finitely generated in every dimension other than four.
The failure of the four-dimensional analogue is the following:

\begin{theo}[Baraglia~\cite{Baraglia23mapping}, Konno~\cite{konno2023homology} \footnote{\cite{Baraglia23mapping} and \cite{konno2023homology} were carried out independently and appeared simultaneously on arXiv.}]
\label{theo: MCG infinite}
There exists a smooth simply-connected closed $4$-manifold $X$ for which the mapping class group $\pi_0(\Diff(X))$ is not finitely generated.
\end{theo}

To give an explicit example of such $X$, let $E(n)$ be a simply-connected elliptic surface over $\CP^{1}$ with $\chi(E(n)) = 12n$ and no multiple fibers.
For $X=E(n)\# S^2\times S^2$ with $n\geq 2$, it was shown in~\cite{Baraglia23mapping,konno2023homology} that $\pi_0(\Diff(X))$ is infinitely generated.
Later, Baraglia--Tomlin~\cite{baraglia2025exoticdiffeomorphisms4manifoldsb} proved that the mapping class group of $X=E(1)\# S^{2}\times S^{2}$ is also infinitely generated.

Let us compare this with the topological category.
First, for a simply-connected closed $4$-manifold $X$, it follows from the results of Freedman~\cite{Fre82} and Quinn~\cite{Q86} that $\pi_0(\Homeo(X))$ is finitely generated (see Remark~\ref{rem Quinn Gabai et al} concerning \cite{Q86}).
Therefore, the infinite generation of the mapping class group of a simply-connected manifold (Theorem~\ref{theo: MCG infinite}) is a phenomenon specific to the smooth category in dimension four.
Note that in the non-simply-connected case, it was already known prior to these results that mapping class groups can be infinitely generated both in higher dimensions and in dimension four (Hatcher~\cite{Hat76}, Budney--Gabai~\cite{budney2021knotted}, Watanabe~\cite{watanabe2023thetagraph}).

In \cite{Baraglia23mapping,konno2023homology}, the infinite generation of the mapping class group is established by showing that its abelianization is infinitely generated, and this approach generalizes to the following statement in higher degrees:

\begin{theo}[Konno \cite{konno2023homology}]
\label{theo: homology infinite}
For every $k\geq 1$, there exists a smooth simply-connected closed $4$-manifold $X$ for which $H_k(\BDiff(X);\Z)$ is not finitely generated.
\end{theo}

More concretely, for $X=E(n)\#_k S^2\times S^2$ with $n\geq 2$, the group $H_k(\BDiff(X);\Z)$ contains $(\Z/2)^{\oplus \N}$ as a direct summand.
Since $H_{1}(\BDiff(X);\Z)$ is the abelianization of $\pi_{0}(\Diff(X))$, Theorem~\ref{theo: homology infinite} with $k=1$ immediately implies Theorem~\ref{theo: MCG infinite}.
Note that for a general simply-connected topological $4$-manifold $X$, finite generation of $\pi_k(\Homeo(X))$ for $k>0$ and of $H_k(\BHomeo(X);\Z)$ for $k>1$ is currently still open.

\subsection{Proof of infinite generation of $H_k(\BDiff(X);\Z)$}

We explain the proof of infinite generation of $H_k(\BDiff(X);\Z)$ (Theorem~\ref{theo: homology infinite}) in somewhat greater detail.
Baraglia's proof of Theorem~\ref{theo: MCG infinite} in \cite{Baraglia23mapping} is also essentially the same as the case $k=1$ of the proof in \cite{konno2023homology}.

Gauge-theoretic characteristic classes are again used in the proof of this infinite generation.
More precisely, one slightly refines the characteristic classes used in the proof of homological instability (Theorem~\ref{theo: homological instability}) to define infinitely many characteristic classes, and then proves that they are linearly independent over $\Z/2$.
To prove this linear independence, one uses the multiple mapping tori that also appeared in the proof of Theorem~\ref{theo: homological instability}.
In our proof, the infiniteness of $H_k(\BDiff(X);\Z)$ ultimately comes from the fact that the elliptic surface $E(n)$ admits infinitely many exotic structures detected by the Seiberg--Witten invariants (see Remark~\ref{rem: Rub infinite}).
Using this infiniteness, one constructs infinitely many multiple mapping tori.
One then evaluates them against the gauge-theoretic characteristic classes to prove linear independence -- this is the overall strategy.
Here, gauge-theoretic characteristic classes capture the ``complexity'' of diffeomorphism groups of $4$-manifolds in the form of infiniteness of the homology of moduli spaces.

Let us describe a little more concretely how these characteristic classes are defined.
The infinitely many characteristic classes are indexed by orbits of the $\Diff^+(X)$-action on $\Spinc(X)/(\Z/2)$, where $\Spinc(X)$ is quotiented by charge conjugation.
More generally, given a $\Diff^+(X)$-invariant subset $\calS\subset \Spinc(X)/(\Z/2)$, one defines a characteristic class
\[
\SWbb^k(X,\calS) \in H^k(\BDiff^+(X);\Z/2).
\]

By collecting the evaluations of these classes, one obtains a homomorphism
\begin{align}
\label{eq: eval infinite}
\bigoplus_{\calS} \langle \SWbb^k(X,\calS), - \rangle : H_{k}(\BDiff^{+}(X);\Z) \to \bigoplus_{\calS}(\Z/2).
\end{align}
In many situations, there are infinitely many choices of $\calS$.
Thus, if one can show that the homomorphism \eqref{eq: eval infinite} surjects onto a subgroup isomorphic to $(\Z/2)^{\oplus \mathbb{N}}$ inside $\bigoplus_{\calS}(\Z/2)$, then it follows that $H_{k}(\BDiff^{+}(X);\Z)$ is not finitely generated.
In the proof of Theorem~\ref{theo: homology infinite} in \cite{konno2023homology}, using computations analogous to those appearing in the proof of Theorem~\ref{theo: homological instability}, it is shown that the homomorphism \eqref{eq: eval infinite} surjects onto the subgroup generated by infinitely many suitably chosen subsets $\calS$.
Such subsets $\calS$ arise from logarithmic transforms of $E(n)$, which produce infinitely many exotic structures on $E(n)$, and infinitely many of them are distinguished by the divisibility of the first Chern class of the spin$^c$ structures, which is an invariant of $\calS$.

\begin{rem}
\label{rem: Rub infinite}
Most of the arguments proving the infiniteness results discussed in this section have as their prototype Ruberman's proof~\cite{Rub99} of infinite generation of the Torelli groups of some $4$-manifolds.
In particular, the idea of constructing families that inherit infinitely many exotic structures of $4$-manifolds and distinguishing them via invariants from gauge theory for families goes back to Ruberman.
With the development of gauge theory for families, new invariants have been introduced, and variants and generalizations of such families have been detected in a variety of different contexts -- including higher homotopy groups of diffeomorphism groups and homology groups of classifying spaces -- leading to the diverse results described above.
Note also that, for the Torelli group, finite generation is known in the higher-dimensional simply-connected case by Sullivan's results~\cite{Sull77}, and hence Ruberman's result~\cite{Rub99} is likewise a phenomenon specific to dimension four.
\end{rem}

\section{Dehn twists and exotic diffeomorphisms}
\label{section exotic Dehn twist}

\subsection{Exotic diffeomorphisms}
\label{subsection exotic diffeo}

Let $X$ be a smooth manifold and $f : X \to X$ a self-diffeomorphism.
If $f$ is topologically isotopic to the identity but not smoothly so, then $f$ is called an {\it exotic diffeomorphism}.
In other words, $f$ is called exotic if it lies in the identity component of $\Homeo(X)$ but not in that of $\Diff(X)$.
Equivalently, $f$ represents a nontrivial element in the kernel of the natural map
\[
\pi_{0}(\Diff(X)) \to \pi_{0}(\Homeo(X)).
\]

In the case with boundary, if $W$ is a manifold with boundary, an element of $\Diff_\del(W)$ that lies in the identity component of $\Homeo_\del(W)$ but not in that of $\Diff_\del(W)$ is called a {\it relatively exotic diffeomorphism}.
Equivalently, it represents a nontrivial element in the kernel of the map
\[
\pi_{0}(\Diff_{\del}(W)) \to \pi_{0}(\Homeo_{\del}(W)).
\]

Exotic diffeomorphisms can be viewed as the most basic manifestation of the discrepancy between the topology of the homeomorphism group and that of the diffeomorphism group.
Such diffeomorphisms have long been known in high dimensions.
For example, given an exotic sphere $\tilde{S}^{n+1}$ in high dimensions, the $h$-cobordism theorem implies that it can be obtained by gluing two copies of $D^{n+1}$ along a boundary diffeomorphism $f : S^{n} \to S^{n}$, and the fact that $\tilde{S}^{n+1}$ is exotic implies that $f$ is exotic.
In dimensions $\leq 3$, the natural map $\Diff(X) \hookrightarrow \Homeo(X)$ is a weak homotopy equivalence.
Therefore, dimension four is the smallest dimension in which exotic diffeomorphisms can occur.

The first example of an exotic diffeomorphism of a $4$-manifold was constructed by Ruberman~\cite{Rub98} in 1998.
This was also the first topological application of gauge theory for families and serves as an important prototype for many of the results discussed in this article.
More recently, four-dimensional analogues of Dehn twists of a quite different nature have been studied actively.
We start by describing these developments.

For the results below, the topological isotopy of the diffeomorphisms to the identity is established by methods independent of gauge theory, based on work of Quinn~\cite{Q86} in the closed case and of Orson--Powell~\cite{OrsonPowell2022} in the case with boundary.

\begin{rem}
\label{rem Quinn Gabai et al}
Recently, a gap in \cite{Q86} was pointed out in \cite{gabai2023pseudoisotopies}, and that paper also filled the gap.
\end{rem}

\subsection{Definition of a Dehn twist}

Let $Y$ be a smooth manifold, and let $\phi : S^1 \to \Diff(Y)$ be a loop based at the identity.
Then $\phi$ defines a diffeomorphism of $Y \times [0,1]$, fixing the boundary pointwise, by
\begin{align}
\label{eq: Dehn twist}
Y \times [0,1] \to Y \times [0,1], \quad (y,t) \mapsto (\phi(t)\cdot y, t).
\end{align}
When $Y = S^1$ and $\phi$ arises from the natural action of $S^1$ on itself, this recovers the classical Dehn twist on $S^1 \times [0,1]$.

Now suppose that $Y$ is a codimension-$1$ submanifold of a manifold $X$ with trivial normal bundle.
By inserting the diffeomorphism \eqref{eq: Dehn twist} into a tubular neighborhood of $Y$, one obtains a diffeomorphism of $X$, called the {\it Dehn twist along $Y$ in $X$}.
If $X$ has boundary $Y$, then after identifying a collar neighborhood of the boundary with $Y \times [0,1]$, the diffeomorphism \eqref{eq: Dehn twist} defines an element of $\Diff_\del(X)$, which is called the {\it boundary Dehn twist} of $X$.

\subsection{Dehn twists along $S^{3}$}

The simplest case is when $Y = S^3$.
In this case, one can consider the Dehn twist associated with a loop representing the generator of
\[
\pi_1(\Diff(S^3)) \cong \pi_1(\SO(4)) \cong \Z/2.
\]
Given $4$-manifolds $X_{1}$ and $X_{2}$, one can perform a Dehn twist along the neck of the connected sum $X_{1}\# X_{2}$, namely, along the cylinder $S^{3} \times [0,1]$ used in forming the connected sum.
As the first example of an ``exotic Dehn twist,'' Kronheimer and Mrowka proved the following:

\begin{theo}[Kronheimer--Mrowka~\cite{KM20Dehn}]
\label{theo: KM Dehn}
The Dehn twist along the neck of the connected sum $K3\# K3$ is an exotic diffeomorphism.
\end{theo}

The proof proceeds by showing that the non-equivariant families Bauer--Furuta invariant is nontrivial for this Dehn twist.
The Bauer--Furuta invariant~\cite{BF04} is a stable homotopy refinement of the Seiberg--Witten invariant, and Kronheimer--Mrowka~\cite{KM20Dehn} explicitly computed its families version in this example.

The idea behind Ruberman's construction of exotic diffeomorphisms~\cite{Rub98} is to construct diffeomorphisms arising from exotic smooth structures on $4$-manifolds.
By contrast, Dehn twists do not arise in any obvious way from exotic smooth structures on $4$-manifolds, and thus represent a qualitatively different source of exotic diffeomorphisms.

Soon afterward, Lin~\cite{Lin20Dehn} proved the following:

\begin{theo}[Lin~\cite{Lin20Dehn}]
\label{theo: Lin Dehn}
The extension of the Dehn twist along the neck of $K3\# K3$ to $K3\# K3\# S^{2}\times S^{2}$ by the identity is also an exotic diffeomorphism.
\end{theo}

Lin's result is proved by computing the $\Pin(2)$-equivariant families Bauer--Furuta invariant of the Dehn twist.
It is known that for a simply-connected $4$-manifold, the mapping class of an exotic diffeomorphism becomes trivial after finitely many stabilizations, that is, after taking connected sums with $S^{2}\times S^{2}$ finitely many times.
In previously known examples, including those of Ruberman~\cite{Rub98} and Baraglia and the author~\cite{BK20gluing}, this triviality occurred after a single stabilization.
The significance of Lin's result is that it provided the first example of an exotic diffeomorphism that survives a single stabilization. More strongly, this is in fact the first example of an exotic phenomenon for simply-connected $4$-manifolds that survives a single stabilization -- including exotic smooth structures on $4$-manifolds themselves.
Later, Kang--Park--Taniguchi~\cite{kang2025exoticdiffeomorphismscontractible4manifold} constructed an example that survives two stabilizations using a Dehn twist along a Seifert fibered $3$-manifold, a construction that we discuss in the next subsection.

Dehn twists along $S^{3}$ have since been studied further; see, for example,
\cite{qiu2024dehntwistconnectedsum,baraglia2024irreducible4manifoldsadmitexotic,tilton2025boundarydehntwistpunctured,lindblad2026boundarydehntwistscommutators}.

\subsection{Dehn twists along Seifert fibered spaces}

To broaden the range of examples, it is natural to consider Seifert fibered spaces as a class of $3$-manifolds.
By definition, such spaces admit an $S^1$-action, and one may consider Dehn twists induced by this action.
It is worth noting that $S^{3}$ is exceptional among Seifert fibered spaces.
Indeed, the loop in $\Diff(S^{3})$ determined by the Seifert $S^{1}$-action is null-homotopic in $\pi_{1}(\Diff(S^{3}))$, and therefore does not define a nontrivial Dehn twist. In other words, the nontrivial loop in $\pi_{1}(\Diff(S^{3}))$ used in the previous subsection does {\it not} come from the Seifert circle action.
By contrast, for a Seifert fibered space $Y\neq S^{3}$, the Seifert $S^1$-action determines a nontrivial element of $\pi_{1}(\Diff(Y))$ (see \cite[Proposition 8.8]{OrsonPowell2022}), and hence potentially gives rise to a nontrivial Dehn twist.

Mallick, Taniguchi, and the author~\cite{konno2023exotic} produced the first examples in which Dehn twists along Seifert fibered spaces give exotic diffeomorphisms.
In the following year, parts of the results of~\cite{konno2023exotic} were significantly generalized in several works: Kang--Park--Taniguchi~\cite{KangParkTaniguchi}, Lin--Mukherjee--Mu\~{n}oz-Ech\'{a}niz and the author~\cite{KLMME2,KLMME}, and Miyazawa~\cite{miyazawaDehn}.
Below, we discuss some results from~\cite{KangParkTaniguchi,KLMME2,KLMME}.

\subsection{Dehn twists on contractible $4$-manifolds.}

First, we summarize results on exotic diffeomorphisms of contractible $4$-manifolds.
The relative topological mapping class group of a compact contractible $4$-manifold is always trivial~\cite{OrsonPowell2022}, and thus, we focus on nontriviality in the smooth category.
The first detection of exotic Dehn twists on contractible $4$-manifolds, carried out in \cite{konno2023exotic}, is proved using a families version of Fr{\o}yshov inequality established in \cite{KT22groupsdiffeomorphismshomeomorphisms4manifolds}.
Another way to detect a Dehn twist is to compute the families Seiberg--Witten invariant of the Dehn twist, which captures phenomena specific to dimension four, described below.
Let $W$ be a smooth $n$-manifold with boundary, and $D^n \hookrightarrow \Int(W)$ be a smooth embedding.
Extending diffeomorphisms by the identity gives a map
\[
i \colon \Diff_\del(D^n) \hookrightarrow \Diff_\del(W).
\]
In high dimensions, if $W$ is contractible, then the diffeomorphism group of $W$ is known to have the same weak homotopy type as that of the disk through the map $i$.
More concretely, the following result holds (dimensions at least six are treated in~\cite{GRW-24-Alexander-trick}, and dimension five in~\cite{krannich2024inftyoperadicfoundationsembeddingcalculus}):

\begin{theo}[Galatius--Randal-Williams~\cite{GRW-24-Alexander-trick}, Krannich--Kupers~\cite{krannich2024inftyoperadicfoundationsembeddingcalculus}]
\label{theo: GRW contractible}
Suppose that $W$ is a compact contractible smooth manifold of dimension $n\geq 5$.
Then, for every smooth embedding $D^n \hookrightarrow \Int(W)$, the map
\[
i \colon \Diff_\del(D^n) \hookrightarrow \Diff_\del(W)
\]
is a weak homotopy equivalence.
\end{theo}

On the other hand, Krushkal--Mukherjee--Powell--Warren~\cite{krushkal2024corksexoticdiffeomorphisms} showed that the four-dimensional analogue of this high-dimensional result (Theorem~\ref{theo: GRW contractible}) fails.
Their proof uses a Ruberman-type exotic diffeomorphism detected in \cite{BK20gluing}, together with a localization theorem established in \cite{krushkal2024corksexoticdiffeomorphisms}.
Here, we say that a diffeomorphism $f\in \Diff_\del(W)$ (or its mapping class $[f]$) \emph{localizes to $D^4$} if there exists an embedding $D^4 \hookrightarrow \Int(W)$ such that the mapping class of $f$ lies in the image of
\[
i_\ast \colon \pi_0(\Diff_\del(D^4)) \to \pi_0(\Diff_\del(W)).
\]
In particular, if there exists a diffeomorphism of $W$ that does not localize to $D^4$, then
$i \colon \Diff_\del(D^4) \hookrightarrow \Diff_\del(W)$
is not a weak homotopy equivalence.
Using a method different from that of~\cite{krushkal2024corksexoticdiffeomorphisms}, one obtains an explicit example of such non-localization:

\begin{theo}[Konno--Lin--Mukherjee--Mu\~{n}oz-Ech\'{a}niz~\cite{KLMME}]
\label{KLMME contractible}
There exists a compact contractible smooth $4$-manifold $W$ with Seifert fibered boundary such that, for the boundary Dehn twist $\tau \in \Diff_\del(W)$, every power $\tau^n$ ($n\neq 0$) defines a nontrivial element of $\pi_0(\Diff_\del(W))$ that does not localize to $D^4$.
\end{theo}

Theorem~\ref{KLMME contractible} shows that the boundary Dehn twist has intrinsic complexity, in the sense that it does not come from $D^{4}$.
The manifold $W$ in Theorem~\ref{KLMME contractible} can be described explicitly as a Mazur manifold, and so Theorem~\ref{KLMME contractible} gives an explicit $4$-manifold and diffeomorphism showing that the four-dimensional analogue of the high-dimensional result (Theorem~\ref{theo: GRW contractible}) fails.

To prove Theorem~\ref{KLMME contractible}, one shows the nontriviality of families Seiberg--Witten invariants for the Dehn twist.
Since the families Seiberg--Witten invariant vanishes for any diffeomorphism that localizes to $D^{4}$, the claimed non-localization follows.
Thus gauge-theoretic invariants once again detect subtle complexities of diffeomorphisms that are specific to dimension four.
A closely related non-localization mechanism also appears in~\cite{krushkal2024corksexoticdiffeomorphisms}.

The first examples of exotic Dehn twists on contractible $4$-manifolds, constructed in~\cite{konno2023exotic}, are given by boundary Dehn twists on contractible $4$-manifolds bounded by Brieskorn spheres $\Sigma(2,3,13)$ and $\Sigma(2,3,25)$.
This result was later extensively generalized as follows:

\begin{theo}[Kang--Park--Taniguchi~\cite{KangParkTaniguchi}]
\label{theo: KPT}
Let $W$ be a compact contractible smooth $4$-manifold whose boundary is a Brieskorn homology sphere $S^{3}$.
Then the boundary Dehn twist $\tau_W$ has infinite order in $\pi_0(\Diff_\partial(W))$.
\end{theo}

As mentioned above, Theorem~\ref{theo: KPT} with the result in \cite{OrsonPowell2022} implies that every power $\tau_W^n$ ($n\neq 0$) is a relatively exotic diffeomorphism.
The proof of Theorem~\ref{theo: KPT} uses a version of equivariant Seiberg--Witten theory due to Baraglia--Hekmati~\cite{BH21,BH-Brieskorn}, together with the technique of {\it homotopy coherent actions} introduced in \cite{KangParkTaniguchi}.
Roughly speaking, this technique produces a family over the classifying space of a cyclic group from the assumption that the boundary Dehn twist is smoothly isotopic to the identity.
They then consider the families relative Bauer--Furuta invariant of this family and derive a contradiction using an equivariant Fr{\o}yshov invariant considered in \cite{BH21,BH-Brieskorn}.

\subsection{Dehn twists on symplectic fillings and monodromy of Milnor fibrations}

The preceding results are topological in nature.
We now turn to results in the symplectic setting.
For symplectic $4$-manifolds with Seifert fibered boundary, we have the following general nontriviality result for boundary Dehn twists:

\begin{theo}[Konno--Lin--Mukherjee--Mu\~{n}oz-Ech\'{a}niz~\cite{KLMME2}]
\label{theo: KLMME 1}
Let $Y$ be a Seifert fibered space with $b_1(Y)=0$, and let $W$ be a compact symplectic filling of $Y$ equipped with the standard contact structure.
If $b^{+}(W)>0$, then the boundary Dehn twist $\tau_W$ is of infinite order in $\pi_0(\Diff_\partial(W))$.
\end{theo}

Again by \cite{OrsonPowell2022}, if $W$ is simply-connected, it follows from Theorem~\ref{theo: KLMME 1} that every power $\tau_W^n$ ($n\neq 0$) is a relatively exotic diffeomorphism.
Here $b^{+}(W)$ is the maximal dimension of a subspace of $H^{2}(W;\R)$ on which the intersection form is positive definite.
The assumption $b^{+}(W)>0$ in Theorem~\ref{theo: KLMME 1} is essential: if one drops it, the conclusion no longer holds.
For example, the Poincar\'{e} homology sphere $\Sigma(2,3,5)$ arises naturally as the boundary of a $4$-manifold $W$ with negative definite intersection form $-E_{8}$ obtained by an $S^{1}$-equivariant plumbing, and the standard $S^{1}$-action on $\Sigma(2,3,5)$ extends over $W$.
From this action, one immediately sees that the boundary Dehn twist of $W$ is trivial in $\pi_0(\Diff_\partial(W))$.

The proof of Theorem~\ref{theo: KLMME 1} proceeds by contradiction.
If one assumes that $\tau_W^n$ is smoothly isotopic to the identity (relative to the boundary), then, following \cite{KangParkTaniguchi}, one constructs a smooth family over the classifying space of a finite cyclic group with fiber $W$.
One then derives a contradiction by analyzing the families relative Bauer--Furuta invariant of this family.
Roughly speaking, the contradiction comes from combining Taubes' nonvanishing theorem for Seiberg--Witten invariants of symplectic $4$-manifolds~\cite{Taubes94sympl} with a vanishing theorem for Floer homology due to Baraglia--Hekmati~\cite{BH-Brieskorn}.
This vanishing result concerns Floer homology of $3$-manifolds obtained as quotients of Seifert fibered spaces by sufficiently large cyclic group actions, based on computations by N\'{e}methi~\cite{Nemethi-05}.
To bridge these nonvanishing and vanishing results and derive a contradiction, one also uses recent developments on the relation between Seiberg--Witten theory and contact structures due to Iida~\cite{Iida21} and Iida--Taniguchi~\cite{Iida-Taniguchi21}.

As a corollary to the above result on symplectic fillings (Theorem~\ref{theo: KLMME 1}), one obtains an interesting consequence concerning the monodromy of Milnor fibrations.
First, we review some basic terminology concerning Milnor fibrations.
Let $f\colon \C^3 \to \C$ be a complex polynomial satisfying $f(0)=0$ and having an isolated singularity at the origin.
Milnor's classical fibration theorem asserts that the restriction
\[
f \colon f^{-1}\bigl(B_\delta(\C)\setminus\{0\}\bigr)\cap B_\epsilon(\C^3)
\to B_\delta(\C)\setminus\{0\}
\quad (1 \gg \epsilon \gg \delta >0)
\]
is a fiber bundle whose boundary family is trivialized.
Here $B_\delta(\C)$ denotes the closed disk in $\C$ of radius $\delta$ centered at the origin.
This bundle is called the {\it Milnor fibration} associated with $f$.
Writing its fiber (the {\it Milnor fiber}) as $M_f$, we denote the monodromy of the fibration by
\[
\mu \in \pi_0(\Diff_\partial(M_f)),
\]
which is one of the most important invariants of the Milnor fibration.
While the action of $\mu$ on homology $H_\ast(M_f)$ has been studied classically, it is also natural to study the mapping class itself.

To study the mapping class of $\mu$, we restrict attention to a slightly smaller class of polynomials.
A polynomial $f$ is called {\it weighted homogeneous} if there exist integers $w_1,w_2,w_3,d>0$ such that
\[
f(t^{w_1}z_1,t^{w_2}z_2,t^{w_3}z_3)=t^d f(z_1,z_2,z_3)
\]
for all $(z_1,z_2,z_3)\in\C^3$ and $t\in\C$.
Important examples of weighted homogeneous polynomials include Brieskorn-type singularities and ADE singularities.
By a Brieskorn-type singularity, we mean a singularity defined by a polynomial of the form
\[
f(z_1,z_2,z_3)=z_1^{p_1}+z_2^{p_2}+z_3^{p_3} \quad (p_i\geq 2).
\]
The ADE singularities---also called rational double points, simple singularities, or du Val singularities---are those given by quotient singularities of the form $\C^2/\Gamma$ for finite subgroups $\Gamma\subset SL(2,\C)$.

For ADE singularities, it is known that the monodromy $\mu$ of the Milnor fibration has finite order in $\pi_0(\Diff_\partial(M_f))$.
This is a consequence of Brieskorn's classical simultaneous resolution theorem~\cite{brieskornADE}.
An important consequence of the above result (Theorem~\ref{theo: KLMME 1}) is that finite-order monodromy occurs only in this situation:

\begin{theo}[Konno--Lin--Mukherjee--Mu\~{n}oz-Ech\'{a}niz~\cite{KLMME2}]
\label{theo: KLMME2}
Suppose that $f\colon \C^3 \to \C$ is a weighted homogeneous polynomial with $f(0)=0$ and an isolated singularity at the origin.
The monodromy $\mu$ of the Milnor fibration of $f$ has finite order in $\pi_0(\Diff_\partial(M_f))$ if and {\bf only if} $f$ defines an ADE singularity.
\end{theo}

The proof is based on the observation that, in the weighted homogeneous case, a suitable power of the monodromy is a boundary Dehn twist, and then reduces the problem to Theorem~\ref{theo: KLMME 1}.
The ADE condition corresponds to the intersection form of $W$ in Theorem~\ref{theo: KLMME 1} being negative definite.

If one considers the topological analogue of Theorem~\ref{theo: KLMME2}, or its higher-dimensional analogue obtained by increasing the number of variables, then both statements are known to fail.
For the topological analogue, see Orson--Powell~\cite{OrsonPowell2022}; in higher dimensions, this follows from results of Krannich--Randal-Williams~\cite{krannich-rw}.
Thus, Theorem~\ref{theo: KLMME2} indicates that even such classical diffeomorphisms as the monodromy of a Milnor fibration exhibit phenomena specific to dimension four and reflect a genuine difference between the smooth and topological categories.

\subsection{Exotic diffeomorphisms of irreducible $4$-manifolds and Donaldson's question}

Above, we discussed exotic diffeomorphisms of $4$-manifolds with boundary, such as Milnor fibers.
By contrast, partly because of the currently available proof techniques, most closed $4$-manifolds on which exotic diffeomorphisms have been detected were obtained as connected sums, for instance manifolds of the form $4\CP^2 \# 21\,\overline{\CP}{}^{2}$ or $K3 \# S^{2}\times S^{2}$~\cite{Rub98,BK20gluing}.
These manifolds admit neither complex nor symplectic structures, and it is not known even whether they carry exotic smooth structures.

By contrast, many central examples of $4$-manifolds in four-dimensional topology are complex surfaces.
After blowdown, such manifolds become {\it irreducible $4$-manifolds}, that is, manifolds admitting no nontrivial connected-sum decomposition.
More precisely, a closed $4$-manifold $X$ is called irreducible if, whenever $X=X_1 \# X_2$, at least one of the $X_i$ is a homotopy $S^4$.
These manifolds serve as the building blocks of four-dimensional topology, and they often carry rich structures such as exotic smooth structures, complex structures, and symplectic structures.

The question of whether irreducible $4$-manifolds can admit exotic diffeomorphisms has long been viewed as a fundamental one in the study of diffeomorphism groups of $4$-manifolds.
It was shown recently that the answer is affirmative:

\begin{theo}[Baraglia--Konno~\cite{baraglia2024irreducible4manifoldsadmitexotic}]
\label{theo: BJK irr}
There exist irreducible closed smooth $4$-manifolds that admit exotic diffeomorphisms.
\end{theo}

More concretely, one can take minimal simply-connected complex surfaces, including various elliptic surfaces and complete intersections.
For example, one may consider elliptic surfaces $E(4n)_{p,q}$ obtained by logarithmic transforms ($n,p,q \geq 1$, except for finitely many pairs $(p,q)$) performed on $E(4n)$, as well as a large class of complete intersections in $\CP^{4}$ cut out by two equations.
The proof uses constraints on smooth $4$-manifolds obtained by Baraglia and the author from families Bauer--Furuta invariants~\cite{BK19}, together with the families index theorem.

The construction of exotic diffeomorphisms in~\cite{baraglia2024irreducible4manifoldsadmitexotic} was not entirely explicit.
Subsequently, Dehn twists along certain Seifert fibered $3$-manifolds in $E(4n)_{p,q}$ were shown to be exotic in~\cite{konno2025constraintslefschetzfibrationsfourdimensional}, thereby providing explicit examples of exotic diffeomorphisms of irreducible $4$-manifolds.
The proof in \cite{konno2025constraintslefschetzfibrationsfourdimensional} uses a generalization to Lefschetz fibrations with $4$-manifold fibers of a constraint on smooth families of $4$-manifolds established in \cite{BK19}.

In fact, this explicit description has an important consequence for a problem in symplectic topology.
A Lagrangian $2$-sphere $S$ in a symplectic $4$-manifold $(X,\omega)$ defines a well-known symplectomorphism $\tau_{S} \in \Symp(X,\omega)$ called the {\it Seidel-Dehn twist} along $S$.
This symplectomorphism $\tau_{S}$ extends the antipodal map on $S$ and is supported near $S$.
The square $\tau_{S}^{2}$ is known to be trivial in the smooth mapping class group $\pi_{0}(\Diff(X))$, but Seidel proved that $\tau_{S}^{2}$ is often nontrivial in the symplectic mapping class group $\pi_{0}(\Symp(X,\omega))$; see Seidel's lecture notes~\cite{Sei08} for details.
In particular, the squared Seidel-Dehn twist $\tau_{S}^{2}$ gives an element of the {\it symplectic Torelli group}
\[
 \pi_{0}(\TSymp(X,\omega))
 = \{[f] \in \pi_{0}(\Symp(X,\omega)) \mid f_{\ast}=\id \text{ on } H_{\ast}(X;\Z)\}.
\]
A well-known question in the field, attributed to Donaldson, asks whether the symplectic Torelli group is generated by squared Seidel-Dehn twists along Lagrangian $2$-spheres when $(X,\omega)$ is a closed simply-connected symplectic $4$-manifold.
It is known that the simple connectivity assumption cannot be dropped~\cite{ArabadjiBaykur25,Smirnov-blowup}, and the answer is affirmative for rational surfaces~\cite{li2022symplectictorelligroupsrational}.
As an important consequence of the detection of exotic diffeomorphisms in \cite{konno2025constraintslefschetzfibrationsfourdimensional}, one obtains a negative answer to Donaldson's question:

\begin{theo}[Konno--Lin--Mukherjee--Mu\~{n}oz-Ech\'{a}niz~\cite{konno2025constraintslefschetzfibrationsfourdimensional}]
\label{theo: KLMME3}
There exist infinitely many simply-connected closed minimal symplectic $4$-manifolds $(X,\omega)$ for which $\pi_{0}(\TSymp(X,\omega))$ is not generated by squared Seidel-Dehn twists.
\end{theo}

Such $X$ is again given by $E(4n)_{p,q}$.
The proof uses the explicit description of exotic diffeomorphisms of $E(4n)_{p,q}$ detected in \cite{konno2025constraintslefschetzfibrationsfourdimensional}: these diffeomorphisms $f$ of $E(4n)_{p,q}$ are symplectomorphisms for suitable symplectic structures $\omega$ on $E(4n)_{p,q}$, and hence define elements of $\pi_{0}(\TSymp(X,\omega))$.
Since squared Seidel-Dehn twists are trivial in $\pi_{0}(\Diff(X))$, the nontriviality of $f$ in $\pi_{0}(\Diff(X))$ immediately implies that $f$ cannot be expressed as a product of squared Seidel-Dehn twists in $\pi_{0}(\TSymp(X,\omega))$.

\begin{rem}
Recently, Du--Li also announced a counterexample to Donaldson's question for the one-point blow-up of a K3 surface.
\end{rem}

\section{Other topics}
\label{section others}

\subsection{Further applications of gauge theory for families}
In this article, we have described only a limited part of the recent developments in applying gauge theory for families to diffeomorphism groups.
There are many other topics that also deserve mention.

\subsubsection{Constraints coming from smooth families of $4$-manifolds}
\label{subsubsec: smooth-family-constraints}

One important topic that could not be discussed sufficiently in this paper is the constraints arising from gauge theory for families of smooth $4$-manifolds.
To recall how classical gauge theory (without families) is applied to $4$-dimensional topology, there are broadly two methodologies:
\begin{enumerate}
\item As explained at the beginning of this article, one defines invariants of $4$-manifolds by counting moduli spaces of solutions to gauge-theoretic partial differential equations, and uses them to distinguish $4$-manifolds.
(More precisely, one may sometimes use finite-dimensional approximations of gauge-theoretic equations, rather than moduli spaces themselves, as invariants.)
Typical examples of this approach include work using the Donaldson invariant~\cite{Do90}, the Seiberg--Witten invariant~\cite{W94}, and the Bauer--Furuta invariant~\cite{BF04}.

\item One studies moduli spaces in order to obtain constraints on classical invariants of $4$-manifolds, most notably on the intersection form.
(More precisely, one may sometimes obtain such constraints from finite-dimensional approximations of gauge-theoretic equations, rather than from moduli spaces themselves.)
Typical examples of this approach include Donaldson's diagonalization theorem~\cite{Do83} and Furuta's $10/8$ inequality~\cite{Fu01}.
\end{enumerate}
What we have explained in most of this article corresponds to the families analogue of (1).
Here we would like to emphasize the families analogue of (2).

The underlying idea is as follows.
Consider a family of gauge-theoretic partial differential equations over a smooth fiber bundle whose fiber is a smooth $4$-manifold $X$.
By analyzing the associated parametrized moduli space, or a family of finite-dimensional approximations of the gauge-theoretic equations, one obtains constraints on classical invariants of the fiber bundle.
Such invariants of a fiber bundle are typically encoded by the action of the monodromy on the intersection form of the fiber.

If one can construct a topological fiber bundle $E \to B$ that fails to satisfy such constraints, then its structure group cannot be reduced from $\Homeo(X)$ to $\Diff(X)$.
That is, if one writes the classifying map of $E$ as $\varphi_{E} : B \to \BHomeo(X)$, then the lifting problem
\begin{align*}
\xymatrix{
     & B\Diff(X)\ar[d]\\
    B \ar@{-->}[ru] \ar[r]_-{\varphi_{E}}  &  B\Homeo(X)
    }
\end{align*}
admits no solution.
Obstruction theory then implies that $\Homeo(X)$ and $\Diff(X)$ differ homotopically.

Along these lines, there are results of Kato--Nakamura--the author~\cite{KKN21}, of Baraglia~\cite{Ba21}, and of Baraglia--the author~\cite{BK19,BKnielsen}; there are also results of Nakamura--the author~\cite{Konno-Nakamura-2023}, giving a $\mathrm{Pin}^{-}(2)$-monopole analogue of \cite{Ba21}, and results of Taniguchi--the author~\cite{KT22groupsdiffeomorphismshomeomorphisms4manifolds} extending \cite{Ba21} to $4$-manifolds with boundary.
In particular, Baraglia's result~\cite{Ba21}, which can be viewed as a families analogue of Donaldson's diagonalization theorem, has broad applicability and has led to many further applications.

A further application of the above idea is the detection of non-smoothable group actions.
In general, when a group $G$ acts on a manifold $X$, the Borel construction yields a family of $X$ over the classifying space $BG$,
\[
EG \times_{G} X \to BG.
\]
Applying constraints from gauge theory for families to this family gives constraints on smooth $G$-actions.
Works that actually detect non-smoothable actions on $4$-manifolds in this way include Nakamura~\cite{Naka10} and Baraglia~\cite{Baraglia19obstructions}.

Some of these themes are explained in more detail in the survey~\cite{KonnoRonsetsu}.
Applications of Baraglia's result~\cite{Ba21} and of Taniguchi--the author's result~\cite{KT22groupsdiffeomorphismshomeomorphisms4manifolds} subsequent to \cite{KonnoRonsetsu} include \cite{konno2023diffeomorphismsexoticphenomenasmall,konno2023exotic,konno2024exoticallyknottedclosedsurfaces,konno2025constraintslefschetzfibrationsfourdimensional}.
These include applications that at first sight have little to do with families, such as the detection of exotic $4$-manifolds and of exotic embeddings of surfaces and $3$-manifolds into $4$-manifolds.

\subsubsection{Applications of a secondary-invariant type}

Another line of applications omitted from this article can be viewed from a perspective reminiscent of secondary invariants.
Namely, when one has a vanishing theorem for solutions of the Seiberg--Witten equations, one can sometimes investigate the nontrivial topology of the ``space of reasons for vanishing.''
Many of the following results are not stated explicitly in that form, but they can be understood from the perspective of such a secondary-invariant viewpoint.

\begin{enumerate}
\item The author's results~\cite{K16,K17} on applications to configurations of surfaces in $4$-manifolds using families of Riemannian metrics.
These are related to the adjunction inequality, an inequality concerning the genus of a surface in a $4$-manifold.
It is known that the Seiberg--Witten equations admit a vanishing theorem for Riemannian metrics obtained by stretching the neighborhood of a surface that violates the adjunction inequality.
(See, for example, \cite{KM94}.)
Using this vanishing theorem, one can impose constraints on configurations of surfaces that violate the adjunction inequality.

\item Applications to surfaces in $4$-manifolds due to Baraglia~\cite{B202}, Iida--Mukherjee--Taniguchi--the author~\cite{iida2022diffeomorphisms}, Auckly~\cite{auckly2023smoothly}, and applications to embeddings of $3$-manifolds in $4$-manifolds due to Iida--Mukherjee--Taniguchi--the author~\cite{iida2022diffeomorphisms} and Mukherjee--Taniguchi--the author~\cite{konno2022exoticcodimension1submanifolds4manifolds}.
These also use vanishing theorems for the Seiberg--Witten equations for metrics obtained by stretching neighborhoods of surfaces violating the adjunction inequality, or neighborhoods of $3$-manifolds with certain special properties.
A typical class of such $3$-manifolds is $L$-spaces.

\item More differential-geometric applications include Ruberman's result~\cite{Rub01} on the topology of the space of positive scalar curvature metrics, the author's result~\cite{K19}, and the result of Auckly--Ruberman~\cite{auckly2025familiesdiffeomorphismsembeddingspositive}.
A vanishing theorem for solutions to the Seiberg--Witten equations for positive scalar curvature metrics has been well known since Witten~\cite{W94}, and this can be used to study the topology of the space of positive scalar curvature metrics.

\item In relation to geometric structures, there are comparisons between symplectomorphism  groups and diffeomorphism groups, and results on the space of symplectic forms due to Kronheimer~\cite{Kropre}, a series of results starting from Smirnov~\cite{Smi20}, and results of Lin~\cite{lin2022family}.
These results ultimately rely on Taubes' vanishing theorem~\cite{Tau95}.
Taubes' nonvanishing theorem for Seiberg--Witten invariants of symplectic $4$-manifolds~\cite{Taubes94sympl} is well known; it guarantees that the Seiberg--Witten invariant is nontrivial for the canonical spin$^{c}$ structure determined by the symplectic structure.
Taubes also showed that, when a spin$^{c}$ structure different from the canonical one satisfies certain conditions, solutions to the Seiberg--Witten equations vanish under Taubes' perturbation using the symplectic form~\cite{Tau95}.
Using this, one can study the topology of the space of symplectic forms.
For further recent developments at the interface of families Seiberg--Witten theory and symplectic geometry, see also~\cite{du2025familyseibergwittenequationkahler,munozechaniz2025configurationslagrangianspheresk3,konno2025constraintslefschetzfibrationsfourdimensional}.
\end{enumerate}

\subsection{Kontsevich characteristic classes}
\label{subsection Kontsevich}
Another major recent direction in the study of diffeomorphism groups of $4$-manifolds, independent of gauge theory for families, stems from Watanabe's decisive work~\cite{Wa18}, which gave a negative answer to the $4$-dimensional Smale conjecture ($\Diff^+(S^4) \not\simeq \mathrm{SO}(5)$).
This relies on Kontsevich characteristic classes, which are quite different in nature from gauge-theoretic invariants, and many further developments have followed.

Among these later developments, we would like to mention the work of Lin--Xie~\cite{lin2023configuration}, which proves that Kontsevich characteristic classes are well defined even for ``formally smooth'' fiber bundles.
In slogan form, their result may be summarized as follows:
``Kontsevich characteristic classes detect aspects of the moduli space of $4$-manifolds that resemble the high-dimensional situation, whereas gauge theory for families detects phenomena specific to dimension four.''
In this sense, \cite{lin2023configuration} suggests that the scope of Kontsevich characteristic classes and that of gauge theory for families are complementary.

To explain this more concretely, we need to introduce some terminology.
Let $X$ be an $n$-dimensional topological manifold, and let $E$ be a topological fiber bundle with fiber $X$.
A {\it formally smooth structure} on $E$ means a reduction of the vertical microtangent bundle $\tau_{v}E$ to $O(n)$, when $\tau_{v}E$ is regarded as a fiber bundle with structure group $\Homeo(\R^{n})$ (by the Kister--Mazur theorem).
Intuitively, this means that, at the level of the vertical tangent bundle, the bundle looks like a smooth fiber bundle.
One can define, in a relatively straightforward way, a classifying space $\calM^{fs}(X)$ for formally smooth fiber bundles; we call it the formally smooth moduli space.
Roughly speaking, $\calM^{fs}(X)$ is defined as
\[
E\Homeo(X) \times_{\Homeo(X)} \{\text{reductions of }\tau X\text{ to }O(n)\}.
\]
Here $\tau X$ denotes the microtangent bundle of $X$.
Of course, any smooth fiber bundle admits a formally smooth structure, and hence there is a forgetful map from the smooth moduli space $\calM^{sm}(X) = \BDiff(X)$,
\[
\calM^{sm}(X) \to \calM^{fs}(X).
\]
Moreover, forgetting the formally smooth structure gives a forgetful map to the topological moduli space $\calM^{top}(X) = \BHomeo(X)$,
\[
\calM^{fs}(X) \to \calM^{top}(X).
\]
Thus the forgetful map $\calM^{sm}(X) \to \calM^{top}(X)$ factors as
\[
\calM^{sm}(X) \to \calM^{fs}(X) \to \calM^{top}(X).
\]
In dimensions other than four, smoothing theory implies that $\calM^{sm}(X) \to \calM^{fs}(X)$ is a weak homotopy equivalence.
Therefore, in dimensions other than four, the entire difference between $\calM^{sm}(X)$ and $\calM^{top}(X)$ comes from the map $\calM^{fs}(X) \to \calM^{top}(X)$.

The result of Lin--Xie~\cite{lin2023configuration} says that, for $X=D^{4}$, an invariant\footnote{More precisely, one considers a relative setting with respect to the boundary.} obtained from Kontsevich characteristic classes, which a priori is defined on the smooth moduli space $\calM^{sm}(X)$, in fact factors through an invariant coming from the formally smooth moduli space $\calM^{fs}(X)$.
It follows that the homotopical difference between $\BDiff(X)=\calM^{sm}(X)$ and $\BHomeo(X)=\calM^{top}(X)$ detected by Watanabe~\cite{Wa18} is actually the homotopical difference of $\calM^{fs}(X) \to \calM^{top}(X)$.
In other words, Watanabe's result detects the difference between $\calM^{sm}(X)$ and $\calM^{top}(X)$ in the same way as in dimensions other than four.

On the other hand, the gauge-theoretic characteristic classes repeatedly discussed in this article are formulated as cohomology classes on $\calM^{sm}(X)$, and one can show that for many concrete examples of $X$ they do not come from $\calM^{fs}(X)$; that is, they do not lie in the image of the natural map $H^{\ast}(\calM^{fs}(X)) \to H^{\ast}(\calM^{sm}(X))$.
As an immediate consequence, the homotopical difference between $\BDiff(X)=\calM^{sm}(X)$ and $\BHomeo(X)=\calM^{top}(X)$ that can be detected by gauge-theoretic characteristic classes is genuinely four-dimensional: it comes from the homotopical difference of $\calM^{sm}(X) \to \calM^{fs}(X)$.

\section{Future directions and outlook}
\label{section future}

There are still many open problems in the gauge-theoretic study of diffeomorphism groups of $4$-manifolds.

For the homology groups of the moduli space $H_\ast(\BDiff(X))$, all classes detected so far by gauge theory for families are torsion.
It is natural to expect that analogous results, such as homological instability and infinite generation, should also hold with $\Q$-coefficients.

Among the directions that remain largely unexplored, no phenomenon specific to dimension four has yet been found at the level of the abstract group structure of $\Diff(X)$, after forgetting its natural topology.
It is an interesting problem whether such phenomena can be discovered by studying the algebraic structure of $\Diff(X)$ through gauge theory.

Moreover, now that exotic diffeomorphisms have been found on irreducible $4$-manifolds (Theorem~\ref{theo: BJK irr}), it may become possible to investigate more deeply the relationships between diffeomorphism groups and various structures on $4$-manifolds, such as exotic smooth structures, complex structures, and symplectic structures.
See~\cite{konno2025constraintslefschetzfibrationsfourdimensional} for an initial attempt in this direction.

In addition, as explained in \ref{subsection Kontsevich}, it seems likely that the complementary study of diffeomorphism groups of $4$-manifolds via Kontsevich characteristic classes and gauge theory for families will continue to develop.

Looking back at the history of the subject in dimensions other than four, the study of diffeomorphism groups of $4$-manifolds will surely remain fascinating for a long time to come.
The subject is still in its infancy.

\bibliographystyle{alpha}
\bibliography{mainref}

\end{document}